\documentclass[reqno,11pt]{amsart}
\usepackage[dvipsnames,usenames]{color} 
\usepackage[toc,page]{appendix}

\usepackage[a4paper,
left=2.5cm, right=2.5cm,
top=3cm, bottom=3cm]{geometry}
    
\usepackage{tikz}  
\usepackage{mathcomp,wasysym}  
\usepackage{mathtools}

\usepackage{hyperref}
\hypersetup{
    colorlinks,
    linkcolor={blue},
    citecolor={blue}
}
\usepackage[mathscr]{euscript}  
\definecolor{luh-dark-blue}{rgb}{0.0, 0.313, 0.608}
         
\usepackage{cite}      
\usepackage{graphicx}  
\usepackage[all]{xy} \xyoption{arc} \xyoption{color}

\usepackage{amsmath} 
\usepackage{amsthm} 
\usepackage{amsfonts}  
\usepackage{amssymb} 
\usepackage{verbatim}  
\usepackage{mathrsfs}

\newcommand{\ZZ}{{\mathbb Z}}
\newcommand{\NN}{{\mathbb N}}
\newcommand{\TT}{{\mathbb T}}
\newcommand{\pat}{\partial_t}
\newcommand{\pax}{\partial_x}

\newcommand{\vertiii}[1]{{\left\vert\kern-0.25ex\left\vert\kern-0.25ex\left\vert #1 
    \right\vert\kern-0.25ex\right\vert\kern-0.25ex\right\vert}}

\newcommand{\bh}{h_\sharp}
\newcommand{\bG}{\Gamma_\sharp}

\newtheorem{lem}{Lemma}

\newtheorem{theorem}{Theorem}  
\newtheorem{definition}{Definition}

\newtheorem{remark}{Remark}

\title[]{On a thin film model with insoluble surfactant}

\author[G. Bruell]{Gabriele Bruell}
\address{Institute for Analysis, Karlsruher Institute of Technology (KIT), D-76128 Karlsruhe, Germany}
\email{gabriele.bruell@kit.edu}

\author[R. Granero-Belinch\'{o}n]{Rafael Granero-Belinch\'{o}n}
\address{Departamento  de  Matem\'aticas,  Estad\'istica  y  Computaci\'on,  Universidad  de Cantabria.  Avda.  Los  Castros  s/n,  Santander,  Spain.}
\email{rafael.granero@unican.es}

\thanks{Date: \today }



\begin{document}

\begin{abstract} This paper studies the existence and asymptotic behavior of global weak solutions for a thin film equation with insoluble surfactant under the influence of gravitational, capillary and van der Waals forces. We prove the existence of global weak solutions for \emph{medium sized} initial data in \emph{large function spaces}. Moreover, exponential decay towards the flat equilibrium state is established, where an estimate on the decay rate can be computed explicitly. 
\end{abstract}

\maketitle 


\section{Introduction}
\allowdisplaybreaks

\emph{Surfactant} is the short form for \emph{surface active agent} and is a substance which -- in contact with a fluid -- reduces surface tension. The induced dynamic is twofold: On the one hand, the resulting surface tension gradients influence the evolution of the thin film; on the other hand, the surfactant speads along the surface. The latter effect is known as \emph{Marangoni effect}. Naturally, the surfactant induced dynamics are of particular interest in connection with \emph{thin fluid films}, where surface tension forces have a very important impact. In particular, the interest in thin film equations with a layer of surfactant on the surface is motived by various applications. For instance coating flow technology, film drainage in emulsions, foams and medical treatment of lungs of premature infants. 

\medskip

The present work studies the dynamics of a viscous, incompressible, Newtonian thin film over a flat bottom equipped with a layer of insoluble surfactant on the free surface. Thus, to study the full problem on has to consider a free boundary problem for the Navier-Stokes equations coupled with an advection-diffusion equation on the free surface. As this is a challenging issue, a common approach to simplify the problem is to consider the \emph{lubrication approximation} to derive evolution equations for the film height and the surfactant concentration which capture the behavior and the main properties of the full free boundary problem. Pioneering works in this direction in absence of surfactant effects are due to Greenspan \cite{greenspan1978motion}, Constantin, Dupont, Goldstein, Kadanoff, Shelley \& Zhou \cite{constantin1993droplet},  Bernis \& Friedman \cite{Bernis1990}, Beretta, Bertsch \& Dal Passo \cite{Beretta1995} and Bertozzi \& Pugh  \cite{Bertozzi1996}. Also, Escher, Matioc \& Matioc \cite{escher2012modelling} considered the flow in porous media (see also Escher \& Matioc \cite{escher2013existence}, Matioc \cite{matioc2012non}, Escher, Lauren\c{c}ot \& Matioc \cite{ELM11}, Lauren\c{c}ot \& Matioc \cite{laurenccot2013gradient,laurencot2014thin,laurenccot2017finite,laurencot2017self}
 and Bruell \& Granero-Belinch\'on \cite{BG18}) while the Stokes flow was considered by Escher, Matioc \& Matioc \cite{escher2013thin} (see also Escher \& Matioc\cite{escher2014non} and Bruell \& Granero-Belinch\'on \cite{BG18}). A more recent reference is Pernas-Casta\~no \& Vel\'azquez \cite{pernas2019analysis}, where the authors study the evolution of the interface between two different fluids in two concentric cylinders when the velocity is given by the Navier-Stokes equation and one of the fluids is thin. 

\medskip
Some of the main works on the evolution of a thin film with insoluble surfactant are the ones by Borgas \& Grotberg \cite{BG}, Gaver \& Grotberg \cite{GG} and Jensen \& Grotberg \cite{JG}. Under certain assumptions, Jensen \& Grotberg \cite{JG} applied the lubrication approximation and cross-sectional averaging to derive the following system of evolution equations for the film height $h=h(t ,x)$ and the surfactant concentration $\Gamma=\Gamma(t,x)$:

\begin{subequations}	\label{eq:system}
	\begin{alignat}{2}
		\partial_t h &=-\partial_x \left[ \frac{h^2}{2}\partial_x \sigma(\Gamma)-\frac{\mathcal{G}}{3}h^3\partial_x h + \frac{\mathcal{S}}{3}h^3\partial_x^3h+\mathcal{A}\frac{\partial_x h}{h} \right]\; && \mbox{in}\quad \Omega_T\\
		\partial_t \Gamma &=- \partial_x \left[\Gamma \left(h \partial_x \sigma(\Gamma)-\frac{\mathcal{G}}{2} h^2 \partial_x h +\frac{\mathcal{S}}{2} h^2 \partial_x^3 h + \frac{3\mathcal{A}}{2} \frac{\partial_x h }{h^{2}}\right)-\mathcal{D}\partial_x \Gamma \right]\;\qquad&&  \mbox{in}\quad \Omega_T.
	\end{alignat}
\end{subequations}

\begin{figure}[h]
	\centering
	\begin{tikzpicture}[domain=0:2*pi, scale=1] 
	\draw[color=black] plot (\x,{0.3*cos(\x r)+1}); 
	\draw[very thick, smooth, variable=\x, luh-dark-blue] plot (\x,{0.3*cos(\x r)+1}); 
	\fill[luh-dark-blue!10] plot[domain=0:2*pi] (\x,0) -- plot[domain=2*pi:0] (\x,{0.3*cos(\x r)+1});
	\draw[very thick,<->] (2*pi+0.4,0) node[right] {$x$} -- (0,0) -- (0,2) node[above] {$z$};
	\draw[very thick,->] (2*pi,0) -- (2*pi,1.3);
	\node[right] at (2*pi,0.5) {$h(t,x)$};
	\coordinate[label=above:{$\Gamma(t,x)$}] (A) at (pi,0.72);
	\fill[color=luh-dark-blue] (A) circle (1.5pt);
	\draw[-] (0,-0.3) -- (0.3, 0);
	\draw[-] (0.5,-0.3) -- +(0.3, 0.3);
	\draw[-] (1,-0.3) -- +(0.3, 0.3);
	\draw[-] (1.5,-0.3) -- +(0.3, 0.3);
	\draw[-] (2,-0.3) -- +(0.3, 0.3);
	\draw[-] (2.5,-0.3) -- +(0.3, 0.3);
	\draw[-] (3,-0.3) -- +(0.3, 0.3);
	\draw[-] (3.5,-0.3) -- +(0.3, 0.3);
	\draw[-] (4,-0.3) -- +(0.3, 0.3);
	\draw[-] (4.5,-0.3) -- +(0.3, 0.3);
	\draw[-] (5,-0.3) -- +(0.3, 0.3);
	\draw[-] (5.5,-0.3) -- +(0.3, 0.3);
	\draw[-] (6,-0.3) -- +(0.3, 0.3);
	\end{tikzpicture} 
	\caption{Scheme of a  thin film flow with insoluble surfactant}\label{Fig1}
\end{figure}
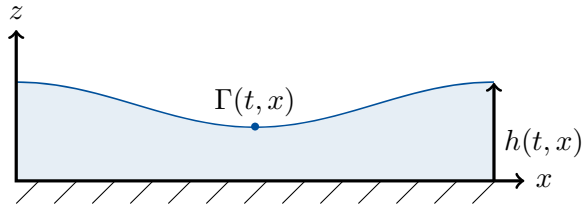
Here, $\Omega_T:=(0,T)\times \Omega$ denotes the time-space domain for the unknown functions $h$ and $\Gamma$, with $\Omega\subset \mathbb{R}$ being an open, bounded interval.
The system \eqref{eq:system} is supplemented with initial conditions
\[
h(0,x)=h_0(x)\quad \mbox{and}\quad \Gamma(0,x)=\Gamma_0(x)\qquad \mbox{for all} \quad x\in \Omega,
\]
where $h_0$ and $\Gamma_0$ are given functions and boundary conditions
$$
	\partial_x h =\partial_x \Gamma=0, \quad \mathcal{S}\partial_x^3 f =0 \qquad \mbox{for all}\quad x\in \partial \Omega.
$$
The appearing parameters represent a modified gravitational constant ($\mathcal{G}$), surface tension coefficient ($\mathcal{S}$), Hamaker constant ($\mathcal{A}$), which corresponds to the effects of van der Waals forces, and surface diffusion coefficient ($\mathcal{D}$). Moreover, $\sigma$ is the constitutive equation of state relating the surface tension to the surfactant concentration. As the presence of surfactant reduces surface tension, $\sigma$ is assumed to be positive and nonincreasing. A commonly used description for the dependence of the surface tension on the surfactant is given by  (cf.  \cite{sheludko1967thin,BG,GG})
\[
\sigma_\beta (s):= (\beta +1) \left[1-s+\left(\frac{\beta +1}{\beta}\right)^\frac{1}{3}s\right]^{-3}-\beta, \qquad s\in[0,1]
\]
for $\beta \in (0,\infty)$. For simplicity reasons, in the present work, we assume that $\sigma$ is given by the limit for $\beta\to  \infty$ in $\sigma_\beta$, that is
$$
\sigma (s)= 1-s, \qquad s\in[0,1].
$$
This assumption is also used in applications and numerical investigations; see for instance \cite{JG2, BGN03, BN04, WCM,HY}.

Let us mention that the scaling Jensen and Grotberg used for the surfactant concentration is given by
\[
\Gamma(t,x)\mapsto \Gamma_m^{-1} \Gamma( t, x),
\]
where $\Gamma_m>0$ is the so-called \emph{critical micelle concentration}. If the surfactant concentration $\Gamma$ exceeds the value $\Gamma_m$, the molecules form micelles and thus there is no further decrease of surface tension to perceive. Consequently, it is natural to assume that the initial surfactant concentration $\Gamma_0$ satisfies $0\leq \Gamma_0\leq 1$.

\medskip

From the analytical point of view, the system \eqref{eq:system} includes many challenges. Notice that the evolution equations in \eqref{eq:system} form a system of two strongly coupled, degenerate, parabolic partial differential equations. Under the assumption that all the appearing parameters $\mathcal G, \mathcal S, \mathcal A,$ and $\mathcal D$ are positive, the degeneracy occurs in the equation for the film height, when $h$ approaches zero, \emph{i.~e.} when the surface touches the bottom. Moreover, \eqref{eq:system} is a coupled system of mixed orders having cross diffusive terms. While the equation for the surfactant concentration $\Gamma$ is an avection-diffusion equation of second order, the equation for the film height $h$ is of fourth order. Notice that, if capillary effects are neglected ($\mathcal{S}=0$), then the system is of second order in both equation. In our considerations, we are going to consider both cases: The gravity driven film ($\mathcal{S}=0$) and the capillary driven film $(\mathcal{S}>0)$.

\medskip

Even if during the last decades modeling as well as numerical investigations for the thin film equation with surfactant have attracted lots of attention (see for instance \cite{GG, JG, JG2, BGN03, BN04, WCM,ECM} and the references therein), the rigorous analytical studies have started recently. Existence of local solutions for a thin film equation with insoluble surfactant driven by Marangoni forces only ($\mathcal G = \mathcal S =\mathcal A=\mathcal D=0$) has been studied by Renardy \cite{R1}. In absence of capillary and van der Waals forces ($\mathcal S=\mathcal A=0$), the authors Escher, Hillairet, Lauren{\c{c}}ot \& Walker \cite{EHLW1} used lubrication approximation to derive a system of differential equations describing the evolution of a thin film with \emph{soluble} surfactant under the influence of Marangoni and gravitational forces. Moreover, they proved local well-posedness in the space of square integrable functions $L^2$ by means of semigroup theory as well as asymptotic stability with exponential decay of the equilibrium. We would like to emphasize that their stability result is stated for positive initial data in the $L^2$-based Sobolev space $H^2$. The result shows in particular, that starting with $H^2$ data close to the flat steady state there exists a unique global strong solution. A similar result for a two-phase thin film equation with insoluble surfactant was shown by Bruell \cite{B1} for the gravity ($\mathcal S =\mathcal A=0$) as well as for the capillary driven film  ($\mathcal G=\mathcal A=0$). 
Due to the degeneracy of the equation with respect to the film height, it is natural to expect that in general strong solutions exist only locally in time. The existence of nonnegative global weak solutions for the thin film equation with insoluble surfactant was proved by Escher, Hillairet, Lauren{\c{c}}ot \& Walker \cite{EHLW11} for the gravity driven thin film  ($\mathcal S=\mathcal A=0$), and \cite{EHLW4th} for the corresponding capillary driven film ($\mathcal G=\mathcal A=0$) and by Bruell \cite{B2} for the two-phase thin film equation with insoluble surfactant under the influence of capillary forces (see also \cite{GW06, CT13}). The main ingredient in all these works concerning the existence of global weak solutions is a regularization argument to overcome the degeneracy, followed by a two-step compactness argument based on a priori estimates provided by an energy functional for the system. Finally, we would like to mention that traveling wave solutions of a gravity thin film equation with insoluble surfactant ($\mathcal S=\mathcal A=0$) were studied by Escher, Hillairet, Lauren{\c{c}}ot \& Walker \cite{EHLWtrav} (see also \cite{EL17}).

\subsection{Aim and outline of the present paper}
The aim of the present work is to prove the existence of global weak solutions of \eqref{eq:system} under fairly low regularity assumptions with respect to the initial data. Similar as in the companion paper \cite{BG18}, we work in scales of Wiener spaces. Exploiting the algebra inequality verified by the norms of the underlying spaces, we show \emph{a priori} energy estimates in Wiener algebra, which guarantee the existence of global weak solutions and imply the exponential decay towards the flat equilibrium state. Moreover, the decay rate can be bounded by explicit constants, which depend on the parameters of the system and the size of the initial data. In addition we prove uniqueness of the weak solutions provided that they belong to a (slightly) higher regularity class. A similar approach has been employed before for the Muskat problem \cite{gan1,gan2,gan3} (and the references therein) for the doubly parabolic Keller-Segel system \cite{burczak2016generalized}, PDEs modelling small steepness porous flow \cite{granero2018asymptotic} and the evolution of crystal surfaces \cite{granero2019global}. We consider two cases: The gravity driven film, where surface tension effects are neglected ($\mathcal S =0$), which leads to a coupled system of second order equations; and the capillary driven film, where we take surface tension effects into account ($\mathcal S>0$). In the latter case the evolution equations \eqref{eq:system} build a coupled system of mixed orders. Let us emphasize that in our work we take \emph{all} acting forces (gravity, surface tension, van der Waals) into account. To the best of our knowledge this is the first analytical existence result for the full system \eqref{eq:system} where $\mathcal G,\mathcal S,\mathcal A,\mathcal D>0$. The outline of the paper is as follows: We start in Section \ref{S:P} with some preliminaries and auxiliary results concerning the scale of Wiener spaces. In Section \ref{S:M} we reformulate the problems in terms of the distance to the (flat) equilibrium and state our main results. Eventually Sections \ref{S:T1} and \ref{S:T2} are devoted to the proofs of the main theorems for the gravity ($\mathcal S=0$) and capillary ($\mathcal S>0$) driven flow, respectively.

\section{Preliminaries}\label{S:P}

We  start by introducing the functional analytical framework. Let $\TT:=[-\pi,\pi)$.
For  $n\in \NN$ we denote by
\begin{equation}\label{Sobinhomo2}
W^{n,p}(\TT)=\left\{u\in L^p(\TT)\text{ such that } \|u\|_{W^{n,p}(\TT)}^p:= \|u\|_{L^p}^p+\|\partial_x^n u\|_{L^p}^p<\infty
\right\}
\end{equation}
the standard $L^p$-based Sobolev spaces on $\TT$. We recall the expression of the $k-$th Fourier coefficient and the Fourier series of a $2\pi-$periodic integrable function $u$,
$$
\hat{u}(k)=\frac{1}{2\pi}\int_{\TT}u(x)e^{-ix k}dx,\;u(x)=\sum_{k\in\ZZ}\hat{u}(k)e^{ix k}.
$$
We introduce the Wiener spaces $\dot{A}^s(\TT)$ as
\begin{equation}\label{Wienerinhomo}
\dot{A}^s(\TT)=\left\{u\in L^1(\TT)\text{ such that } \|u\|_{\dot{A}^s(\TT)}:=\sum_{k\in\ZZ} |k|^s|\hat{u}(k)| < \infty\right\}.
\end{equation}
We note that $\dot{A}^0(\TT)=A(\TT)$ is a Banach algebra and $\{\dot{A}^s(\TT)\mid s\geq 0\}$ form a Banach scale. Furthermore, the following inequalities hold true:
\begin{itemize}
\item	Let $p\geq q\geq 0$, then
\begin{equation}
	\label{lem:P}
	\|f\|_{ \dot A^q }\leq \|f\|_{ \dot A^p} \text{ for all } f\in  \dot A^p(\TT).
\end{equation}
\item Let $s\in \{0\}\cup[1,\infty)$ be a fixed parameter and $f,g\in \dot A^s(\TT)$, then, due to the convexity of $x^s$ in this range,
\begin{equation}\label{eq:L21}
\|fg\|_{\dot{A}^s}\leq 2^{s-1}\left\|f\|_{\dot A^s}\|g\|_{\dot A^0}+\|f\|_{\dot A^0}\|g\|_{\dot A^s}\right)\leq  2^{s}\|f\|_{\dot{A}^s}\|g\|_{\dot{A}^s},
\end{equation}
while, as a consequence of subadditivity of $x^s$, for $s\in (0,1)$, we have 
\[
	\|fg\|_{\dot{A}^s}\leq 2\|f\|_{\dot A^s}\|g\|_{\dot A^s}.
\]
\item Due to the H\"older inequality, we have the following interpolation inequality
\begin{equation}\label{interpolation}
\|f\|_{\dot{A}^{s\theta}}\leq \|f\|_{\dot A^0}^{1-\theta}\|f\|_{\dot{A}^{s}}^{\theta}\quad \mbox{for all }\quad 0<\theta<1.
\end{equation}
\end{itemize}

\section{Reformulation of the problem and main results}\label{S:M}
\subsection{Reformulation}
In what follows we assume that $\sigma(s)=1-s$. Then, system \eqref{eq:system} is given by 
\begin{subequations}\label{eq:systemA}
	\begin{alignat}{2}
	\partial_t h &=-\partial_x \left[ -\frac{h^2}{2}\partial_x\Gamma-\frac{\mathcal{G}}{3}h^3\partial_x h + \frac{\mathcal{S}}{3}h^3\partial_x^3h+\mathcal{A}\frac{\partial_x h}{h} \right],\; &&\mbox{in}\quad \Omega_T \\
	\partial_t \Gamma &=- \partial_x \left[\Gamma \left(-h \partial_x \Gamma-\frac{\mathcal{G}}{2} h^2 \partial_x h +\frac{\mathcal{S}}{2} h^2 \partial_x^3 h + \frac{3\mathcal{A}}{2} \frac{\partial_x h }{h^{2}}\right)-\mathcal{D}\partial_x \Gamma \right],\;&&\mbox{in}\quad \Omega_T,
	\end{alignat}
\end{subequations}
with initial conditions
\[
h(0,x)=h_0(x)\quad \mbox{and}\quad \Gamma(0,x)=\Gamma_0(x)\qquad \mbox{for all} \quad x\in \Omega,
\]
and  boundary conditions
\begin{equation}\label{eq:boundary}
\partial_x h =\partial_x \Gamma=0, \quad \mathcal{S}\partial_x^3 f =0 \qquad \mbox{for all}\quad x\in \partial \Omega.
\end{equation}
 The problem is posed on a spatial interval $\Omega=(0,L)$ with the above lateral boundary conditions. However, without lossing generality, instead of considering an interval and no-flux boundary conditions, we are going to consider periodic solutions  $h,\Gamma$ of \eqref{eq:systemA} on a flat torus $\TT$ (which can be identified with $[-\pi,\pi)$). This generalization actually simplifies our approach and it was already used in \cite{BG18} for similar problems. Let us explain why our formulation in the flat torus is actually equivalent to the original problem posed on the interval $(0,L)$. If $(h_0,\Gamma_0)$ are the initial data on an interval $\overline \Omega=[0,L]$ satisfying the boundary conditions \eqref{eq:boundary}, we set
	\[
	\bar h_0(x):= h_0(|x|) ,\quad 	\bar \Gamma_0(x):=\Gamma_0(|x|)\qquad \mbox{for}\quad x\in [-L,L].
	\]
		In view of the symmetry of \eqref{eq:systemA}, the evenness of initial data is preserved and any solution of \eqref{eq:systemA} on $[0,L]$ with initial data $(h_0,\Gamma_0)$ satisfying the boundary conditions \eqref{eq:boundary}, can be identified with the corresponding solution to even initial data $(\bar h_0, \bar \Gamma_0)$ on the periodic cell $[-L,L]$ restricted to the half-domain $[0,L]$. In the sequel we drop the bar notation and consider periodic solutions of \eqref{eq:systemA} defined on $\TT$ with initial conditions 
\[
h(0,x)=h_0(x)\quad \mbox{and}\quad \Gamma(0,x)=\Gamma_0(x)\qquad \mbox{for all} \quad x\in \TT,
\]
where $h_0,\Gamma_0$ are given periodic functions.
It follows immediately from the structure of the equations \eqref{eq:systemA} that the initial mass is preserved in time:

\begin{lem}[Conservation of mass]
Let $(h,\Gamma)$ be a solution of \eqref{eq:systemA} on a time interval $[0,T)$, then
\[
	\int_{\TT} h(t,x)\,dx =	\int_{\TT} h_0(x)\,dx\quad\mbox{and}\quad \int_{\TT} \Gamma(t,x)\,dx =	\int_{\TT} \Gamma_0(x)\,dx\qquad \mbox{for all}\quad t\in[0,T).
\]
\end{lem}

If $(h_0,\Gamma_0)$ are nonnegative bounded initial data, we set 
 $$
 \bh=\frac{1}{2\pi}\int_\TT h_0(x)dx\quad\mbox{and}\quad \bG=\frac{1}{2\pi}\int_\TT \Gamma_0(x)dx.
 $$
 The constants $\bh$ and $\bG$ represent the mean of the initial data and they are a steady state of the system \eqref{eq:systemA}. In our studies we consider the evolution of the distance of a solution $(h,\Gamma)$ to the steady state $(\bh,\bG)$. For this purpose, we define new unknowns
\begin{equation}\label{newunknown}
f=h-\bh,\quad \Theta=\Gamma-\bG,
\end{equation}
which have zero mean. 
In the new variables \eqref{newunknown}, the system \eqref{eq:systemA} can be rewritten as
\begin{subequations}\label{eq:systemlinearized}
\begin{alignat}{2}
\partial_t f-\frac{\bh^2}{2}\partial_x^2 \Theta +\left(\frac{\mathcal{A}}{\bh}-\frac{\mathcal{G}}{3}\bh^3\right)\partial_x^2 f+\frac{\mathcal{S}}{3}\bh^3\partial_x^4 f&=\sum_{j=1}^4 N_j,\; &&\mbox{in}\quad (0,T)\times \TT\\
\partial_t \Theta -\left(\bh\bG+\mathcal{D}\right)\pax^2\Theta+\left(\frac{3\mathcal{A}\bG}{2\bh^2}-\frac{\mathcal{G}}{2}\bG\bh^2\right) \pax^2f+\frac{\mathcal{S}}{2}\bG\bh^2 \pax^4f &= \sum_{j=5}^8 N_j,\;&& \mbox{in}\quad (0,T)\times \TT
\end{alignat}
\end{subequations}
with initial conditions
\[
f(0,x)=f_0(x)\quad \mbox{and}\quad \Theta(0,x)=\Theta_0(x)\qquad \mbox{for all} \quad x\in \Omega,
\]
where $f_0(x)=h_0(x)-\bh$ and $\Theta_0(x)=\Gamma_0(x)-\bG$ are the initial displacement functions from the flat states $\bh$ and $\bG$, respectively.
The nonlinear terms $N_i$, $i=1,\dots,8,$ on the right hand side of \eqref{eq:systemlinearized} are defined as
\begin{align*}
N_1&=\pax\left[\left(\frac{f^2}{2}+f\bh\right)\partial_x \Theta\right],\\
N_2&=\partial_x \left[\frac{\mathcal{G}}{3}\left(3\bh^2 f+3f^2\bh+f^3\right)\partial_x f\right],\\
N_3&=-\partial_x \left[\frac{\mathcal{S}}{3}\left(3\bh^2 f+3f^2\bh+f^3\right)\partial_x^3 f\right],\\
N_4&=\partial_x \left[\mathcal{A}\frac{f}{\bh^2(1+\frac{f}{\bh})}\partial_x f \right],\\
N_5&=\pax\left[\left(\bG f+\Theta\bh+\Theta f\right)\partial_x \Theta\right],\\
N_6&=\pax\left[\frac{\mathcal{G}}{2}\left(\bG f^2+2\bG\bh f+\Theta\bh^2+\Theta f^2+2\Theta\bh f\right)\partial_x f\right],\\
N_7&=-\pax\left[\frac{\mathcal{S}}{2}\left(\bG f^2+2\bG\bh f+\Theta\bh^2+\Theta f^2+2\Theta\bh f\right)\partial_x^3 f\right],\\
N_8&=\partial_x \left[\Theta \frac{3\mathcal{A}}{2\bh^2} \frac{(f^2+2\bh f)\partial_x f }{(\bh+f)^{2}}+\bG\frac{3\mathcal{A}}{2\bh^2} \frac{(f^2+2\bh f)\partial_x f }{(\bh+f)^{2}}-\Theta \frac{3\mathcal{A}}{2\bh^2} \partial_x f \right]\\
&= \partial_x \left[ \frac{3\mathcal{A}}{2}\frac{2f\bG\bh+f^2\bG-\Theta \bh^2}{\bh^2(f+\bh)^2}\partial_x f\right].
\end{align*}

Note that, for $|r|<1$, we have
\begin{align}
\frac{1}{1+r}&=\sum_{j=1}^\infty (-1)^{j+1}r^{j-1},\label{eq:taylor1} \\
\frac{1}{\left(1+r\right)^2}&=\sum_{j=1}^\infty j(-1)^{j+1}r^{j-1},\label{eq:taylor2}\\
\frac{1}{\left(1+r\right)^3}&=\frac{1}{2}\sum_{j=2}^\infty j(j-1)(-1)^{j}r^{j-2}.\label{eq:taylor3}
\end{align}
Consequently, under the assumption that
\begin{equation}\label{restriction1}
\|f\|_{L^\infty}\leq \|f\|_{\dot A^0}< \bh,
\end{equation}
 we can use \eqref{eq:taylor1} and \eqref{eq:taylor2} and write
\begin{align}
N_4&=\mathcal{A} \bigg{[}\left(\frac{f}{\bh^2}\partial_x^2 f +\left(\frac{\partial_x f}{\bh}\right)^2\right)\sum_{j=1}^\infty (-1)^{j+1}\left(\frac{f}{\bh}\right)^{j-1}-\frac{f}{\bh^3}(\partial_x f)^2\sum_{j=1}^\infty j(-1)^{j+1}\left(\frac{f}{\bh}\right)^{j-1} \bigg{]}.\label{NL4}
\end{align}
Similarly, invoking  \eqref{eq:taylor2} and \eqref{eq:taylor3} we have that
\begin{align}
N_8&=\partial_x \left[ \frac{3\mathcal{A}}{2}\frac{2f\bG\bh+f^2\bG-\Theta \bh^2}{\bh^2(f+\bh)^2}\partial_x f\right]\nonumber\\
&= \frac{3\mathcal{A}}{2\bh^4}\left[\frac{2f\bG\bh+f^2\bG-\Theta \bh^2}{\left(\frac{f}{\bh}+1\right)^2}\partial_x^2 f+\frac{2\pax f\bG\bh+2f\pax f\bG-\pax \Theta \bh^2}{\left(\frac{f}{\bh}+1\right)^2}\partial_x f\right]\nonumber\\
&\quad- \frac{3\mathcal{A}}{\bh^4}\left[\frac{2f\bG\bh+f^2\bG-\Theta \bh^2}{\left(\frac{f}{\bh}+1\right)^3}\frac{(\partial_x f)^2}{\bh}\right]\nonumber\\
&=\frac{3\mathcal{A}}{2\bh^4}\left[\left(2f\bG\bh+f^2\bG-\Theta \bh^2\right)\partial_x^2 f+\left(2\pax f\bG\bh+2f\pax f\bG-\pax \Theta \bh^2\right)\partial_x f\right]\sum_{j=1}^\infty j(-1)^{j+1}\left(\frac{f}{\bh}\right)^{j-1}\nonumber\\
&\quad- \frac{3\mathcal{A}}{\bh^4}\left[\left(2f\bG\bh+f^2\bG-\Theta \bh^2\right)\frac{(\partial_x f)^2}{\bh}\right]\frac{1}{2}\sum_{j=2}^\infty j(j-1)(-1)^{j}\left(\frac{f}{\bh}\right)^{j-2}.\label{NL8}
\end{align}

\medskip
We fix the initial data $h_0$, $\Gamma_0$ for problem \eqref{eq:systemA}. Thereby, the constants $\bh$ and $\bG$ are uniquely determined and we are going to state our results in terms of  $f$,$\Theta$ for \eqref{eq:systemlinearized}. 

\subsection{Main results}
In what follows the constants $\mathcal G,\mathcal A,$ representing the gravitational and van der Waals forces as well as the diffusion coefficient $\mathcal D$ are assumed to be strictly positive.

\begin{definition}\label{defi1}Set $\zeta=3$ for the capillary driven flow ($\mathcal{S}>0$) and $\zeta=1$ for the gravity driven flow ($\mathcal{S}=0$).We say that $(f,\Theta)\in \left(L^1(0,T;W^{\zeta,1}(\TT))\right)^2$ is a weak solution of (\ref{eq:systemlinearized}) on $[0,T)$ corresponding to initial data $(f_0,\Theta_0)$ if and only if
\begin{align*}
-\int_\TT f_0\phi(0)dx&-\int_0^T\int_\TT f\pat\phi dxdt+\int_0^T\int_\TT \pax^2\phi\left(-\frac{\bh^2}{2}\Theta +\left(\frac{\mathcal{A}}{\bh}-\frac{\mathcal{G}}{3}\bh^3\right)f\right)  dxdt\\
&+\int_0^T\int_\TT\pax^4\phi \frac{\mathcal{S}}{3}\bh^3f dxdt+\sum_{j=1}^4\int_0^T\int_\TT  \partial_x^{-1}N_j\pax \phi dxdt =0
\end{align*}
and
\begin{align*}
-\int_\TT \Theta_0\psi(0)dx&-\int_0^T\int_\TT \Theta\pat\psi dxdt+\int_0^T\int_\TT \pax^2\psi\left(-\left(\bh\bG+\mathcal{D}\right)\Theta+\left(\frac{3\mathcal{A}\bG}{2\bh^2}-\frac{\mathcal{G}}{2}\bG\bh^2\right)f\right)  dxdt\\
&+\int_0^T\int_\TT\pax^4\psi\frac{\mathcal{S}}{2} \bG\bh^2f dxdt+\sum_{j=5}^8\int_0^T\int_\TT  \partial_x^{-1}N_j\pax \psi dxdt =0
\end{align*}
for all $(\phi,\psi)\in C^{\infty}_c([0,T)\times \TT)$,  where $\pax^{-1}$ denotes the operator given by $\widehat{\pax^{-1}u}(n)=-\frac{i}{n}\hat{u}(n)$. 

\end{definition}

Let $r,s\geq 0$. For $(f,\Theta):[0,T) \to \dot A^r(\TT)\times \dot A^s(\TT)$, we define the following functional:
\begin{equation}\label{energy}
\mathscr{E}^r_{s}(f,\Theta):[0,T) \to [0,\infty),\qquad \mathscr{E}^r_{s}(f,\Theta)(t):=\|f(t)\|_{\dot{A}^r}+\|\Theta(t)\|_{\dot{A}^s}.
\end{equation}

\medskip

We start by formulating our main result for the gravity driven flow, that is when $\mathcal S=0$. To this end, let us define
\begin{align}
\Lambda_1(t)&:= \bh+\frac{19\bh^2\mathcal{ G}}{3}+ \frac{\mathcal{A}}{\bh^2}\left(1-\frac{\mathscr{E}_0^0(f,\Theta)(t)}{\bh}\right)^{-1}\bigg{[}2+\left(1-\frac{\mathscr{E}_0^0(f,\Theta)(t)}{\bh}\right)^{-1} \bigg{]}+\bG\nonumber \\ 
&\quad+\mathcal{G}\left(4\bG\bh +5\bh^2\right)+\frac{3\mathcal{A}}{2\bh^3}\left(1-\frac{\mathscr{E}_0^0(f,\Theta)(t)}{\bh}\right)^{-2}\Big\{7\bG + \frac{3}{2}\bh + \left(1-\frac{\mathscr{E}_0^0(f,\Theta)(t)}{\bh}\right)^{-1}4\bG\Big\} \label{eq:Lambda_1},\\
\Lambda_2(t)&:=\frac{13}{2}\bh+2\bG+\mathcal{G}\bh^2+\frac{3\mathcal{A}}{2\bh^3}\left(1-\frac{\mathscr{E}_0^0(f,\Theta)(t)}{\bh}\right)^{-2}\frac{\bh}{2}\label{eq:Lambda_2}\,,
\end{align}
and the constants
\begin{align}
\mathfrak{C}_1&:=\frac{\mathcal{G}}{3}\bh^3-\frac{\mathcal{A}}{\bh}-\left|\frac{3\mathcal{A}\bG}{2\bh^2}-\bG\bh^2\frac{\mathcal{G}}{2}\right|,\label{constant1}\\
\mathfrak{C}_2&:=\bh \bG+\mathcal{D}-\frac{\bh^2}{2}.\label{constant2}
\end{align}

Then, we have the following result:
\begin{theorem}[Global existence for $\mathcal{S}=0$]\label{global1}
Let $f_0\in \dot A^0(\TT), \Theta_0\in \dot A^0(\TT)$ be nontrivial  initial data for \eqref{eq:systemlinearized}
such that
$$
\mathscr{E}^0_0(f_0,\Theta_0)<\min\{\bh,\bG\}.
$$
Assume that $\mathfrak{C}_1,\mathfrak{C}_2>0$. If
\begin{align*}
\gamma_1:=\mathfrak{C}_1-\Lambda_1(0)\mathscr{E}_0^0(f_0,\Theta_0)&>0,\\
\gamma_2:=\mathfrak{C}_2-\Lambda_2(0)\mathscr{E}_0^0(f_0,\Theta_0)&>0,
\end{align*}
then there exists at least one global solution $(f,\Theta)$ of \eqref{eq:systemlinearized} in the sense of Definition \ref{defi1} with regularity
\begin{align*}
(f,\Theta)\in &\Big(L^{\frac{2}{r}}\left(0,T;W^{r,\infty}(\TT)\right)\cap L^1\left(0,T;C^{1+\alpha}(\TT)\right) \cap L^2(0,T;H^2(\TT))\Big)^2,
\end{align*}
for any $r\in[0,2), \alpha\in[0,\frac{1}{2})$ and $T>0$. Moreover, the weak solution satisfies the following exponential decay:
$$
\|f(t)\|_{L^\infty}+\|\Theta(t)\|_{L^\infty}\leq \mathscr{E}^0_0(f_0,\Theta_0)e^{-\delta t},
$$
where 
\begin{equation}\label{delta2}
\delta=\min\{\gamma_1,\gamma_2\}>0.
\end{equation}
Furthermore, if the solution $(f,\Theta)$ satisfies the regularity
$$
(f,\Theta)\in \left(L^1\left(0,T;\dot{A}^2(\TT)\right)\right)^2,
$$
then the weak solution is unique.
\end{theorem}

\medskip

\begin{remark}
The condition $
\mathscr{E}^0_0(f_0,\Theta_0)<\min\{\bh,\bG\}
$ in Theorem \ref{global1} implies in particular that
\[
	\|f\|_{L^\infty}<\|f_0\|_{\dot A^0}<\bh\quad\mbox{and}\quad \|\Theta\|_{L^\infty}<\|\Theta_0\|_{\dot A^0}<\bG.
\]
This corresponds to a positivity condition of $h_0$ and $\Gamma_0$ and ensures that \eqref{restriction1} is initially satisfied.
\end{remark}

\begin{remark}
Note that the size restriction is explicitly computable and that the initial data do not need to be small in Sobolev spaces (one should compare the results in this paper with the global result in \cite{bruellwell}). In particular, we prove the existence global weak solutions and their decay towards the flat state for highly oscillating initial data. Let us further explain this with an explicit example. Consider the case where van der Waals forces are neglected, that is $\mathcal A=0$ and set $\mathcal G=\mathcal D=1$. We choose the initial data $(h_0,\Gamma_0)$ to be
$$
h_0(x)=1+\mu\sin(1000x)\quad \mbox{and}\quad \Gamma_0(x)=\frac{1}{2}+\mu\cos(1000x),
$$
for $0<\mu< \frac{1}{4}$. Then, 
$$
\bh=1,\quad\bG=\frac{1}{2},\quad\mathscr{E}^0_0(h_0-\bh,\Gamma_0-\bG)=2\mu,
$$
and the constants $\mathfrak{C}_1,\mathfrak{C}_2,\Lambda_1(0),$ and $\Lambda_2(0)$ are given by
$$
\mathfrak{C}_1=\frac{1}{12},\quad\mathfrak{C}_2=1,\quad \Lambda_1(0)= \frac{89}{6},\quad \Lambda_2(0)=\frac{17}{2}.
$$
Then, for $\mu<\frac{1}{356}\sim 0.003$, this family of initial data satisfies the hypotheses and Theorem~\ref{global1} guarantees the existence of a global weak solution $(h,\Gamma)$ corresponding to the initial data $(h_0,\Gamma_0)$. Moreover the solution decays exponentially towards the flat equilibrium and
\[
\|h(t)-\bh\|_{L^\infty}+\|\Gamma-\bG\|_{L^\infty}\leq 2\mu e^{-\delta t},
\]
where $\delta =\frac{1}{12}-\frac{89}{3}\mu$. Furthermore, observe that
$$
\|h_0\|_{\dot{C}^1}=\|\Gamma_0\|_{\dot{C}^1}=O(1) \quad \mbox{while}\quad \|h_0\|_{H^2}=\|\Gamma_0\|_{H^2}=O(10^3).
$$
\end{remark}

\medskip

Next we formulate our main theorem for the capillary driven flow, that is for the full system \eqref{eq:systemA} with $\mathcal S>0$.
In addition to \eqref{eq:Lambda_1}--\eqref{constant2}, we define
{\begin{align}
\Lambda_3&:=\frac{\mathcal{S}}{2} (14\bh\bG +4\bh^2)+\frac{19}{3}\mathcal{S}\bh^2,\label{eq:Lambda_3}\\
\mathfrak{C}_3&:=\frac{\mathcal{S}}{3}\bh^3-\bG\bh^2\frac{\mathcal{S}}{2}\label{constant3}.
\end{align}

Then, we have the following result:
\begin{theorem}[Global existence for $\mathcal{S}>0$]\label{global2} 
Let $f_0\in \dot A^0(\TT), \Theta_0\in \dot A^0(\TT)$ be nontrivial initial data for \eqref{eq:systemlinearized}
such that
$$
\mathscr{E}^0_0(f_0,\Theta_0)<\min\{\bh,\bG\}.
$$
Assume that $\mathfrak{C}_1,\mathfrak{C}_2,\mathfrak{C}_3>0$. If
\begin{align*}
\gamma_1:=\mathfrak{C}_1-\Lambda_1(0)\mathscr{E}_0^0(f_0,\Theta_0)&>0,\\
\gamma_2:=\mathfrak{C}_2-\left(\Lambda_2(0)+\frac{\mathcal{S}\bh^2}{2}\right)\mathscr{E}_0^0(f_0,\Theta_0)&>0,\\
\gamma_3:=\mathfrak{C}_3-\Lambda_3\mathscr{E}_0^0(f_0,\Theta_0)&>0,
\end{align*}
then there exists at least one global solution $(f,\Theta)$ of \eqref{eq:systemlinearized} in the sense of Definition \ref{defi1} with regularity
\begin{align*}
f\in & L^\frac{4}{s}\left(0,T; W^{s,\infty}(\TT)\right)\cap L^{1}\left(0,T;C^{3+\alpha}(\TT)\right) \cap L^2\left(0,T;H^2(\TT)\right),\quad s\in[0,4),\\
\Theta\in & L^\frac{4}{r}\left(0,T; W^{r,\infty}(\TT)\right)\cap L^{1}\left(0,T;C^{1+\alpha}(\TT)\right)\cap L^2\left(0,T;H^1(\TT)\right),\quad r\in[0,2)
\end{align*}
for any $T$ and $0\leq\alpha<\frac{1}{2}$. Moreover, the weak solutions satisfies the following exponential decay:
$$
\|f(t)\|_{L^\infty}+\|\Theta(t)\|_{L^\infty}\leq \mathscr{E}^0_0(f_0,\Theta_0)e^{-\delta t},
$$
where 
\begin{equation}\label{delta}
\delta=\min\{\gamma_1,\gamma_2,\gamma_3\}>0.
\end{equation}
Furthermore, if the solution $(f,\Theta)$ satisfies the regularity
$$
(f,\Theta)\in L^1\left(0,T;\dot{A}^4(\TT)\right)\times L^1\left(0,T;\dot{A}^2(\TT)\right),
$$
then the weak solution is unique.
\end{theorem}

The remainder or this paper is devoted to the proof of the above theorems.

\section{Proof of Theorem \ref{global1}: Global existence and decay when $\mathcal{S}=0$} 
\label{S:T1}
Let us start by outlining the steps of the proof. First, we provide in Subsection \ref{ss:e} some \emph{a priori} estimates of a solution $(f,\Theta)\in \left(C^1(0,T;\dot A^0(\TT))\right)^2$ of \eqref{eq:systemlinearized}. Under the assumptions on the initial data in Theorem \ref{global1}, which in particular require that
\[
	\mathscr{E}_0^0(f_0,\Theta_0)<\min\{\bh,\bG\},
\] 
we will show that the solution $(f,\Theta)$ preserves this estimate, that is
$$
\mathscr{E}_0^0(f,\Theta)(t)<\min\{\bh,\bG\}\quad\mbox{for all}\quad 0<t<T.
$$
Notice that the above inequality  implies that
$$
\|f(t)\|_{\dot A^0}\leq \mathscr{E}_0^0(f,\Theta)(t)<\min\{\bh,\bG\}<\bh \quad\mbox{for all}\quad 0<t<T,
$$
$$
\|\Theta(t)\|_{\dot A^0}\leq \mathscr{E}_0^0(f,\Theta)(t)<\min\{\bh,\bG\}<\bG \quad\mbox{for all}\quad 0<t<T.
$$
Hence, the positivity conditions for $h$ and $\Gamma$, 
$$
\|f(t)\|_{L^\infty}\leq \|f(t)\|_{\dot A^0} < \bh,\quad \|\Theta(t)\|_{L^\infty}\leq \|\Theta(t)\|_{\dot A^0}<\bG 
$$
are preserved for all $t\in [0,T)$. Then we obtain that, for small enough initial energy $\mathscr{E}^0_0(f_0,\Theta_0)$, the following inequality holds:
$$
\frac{d}{dt} \mathscr{E}^0_0(f,\Theta)(t)+\delta \mathscr{E}^2_2(f,\Theta)(t)\leq 0
$$
for some $\delta>0$ dependent on the initial data.
This inequality implies that a local solution $(f,\Theta)$ can not leave a ball in $(L^\infty((0,T); \dot A^0(\TT))\cap L^1(0,T;\dot A^2(\TT)))^2$; thus, there is no finite time singularity in these functional spaces. 
In Subsection \ref{ss:gl} we implement a standard Galerkin approximation argument to prove the existence of global weak solutions. The existence of a classical local solution of the Galerkin approximated sytem is guaranteed by the Picard--Lindel\"of Theorem.  These approximated solutions $(f_M,\Theta_M)_{M\in \NN}$ are analytic. In particular, they satisfy 
$$
(f_M,\Theta_M)\in \left(C^1([0,T_M); \dot A^0(\TT))\right)^2,
$$ where $T_M>0$ denotes the maximal time of existence of $(f_M,\Theta_M)$. The \emph{a priori} estimates from before ensure the \emph{global} existence of the approximated solutions and provide the necessary \emph{a priori} bounds in order to use compactness arguments and pass to the limit $M\to \infty$, which yields a global weak solution of \eqref{eq:systemlinearized} in the sense of Definition \ref{defi1}. 
Now, let $(f_M,\Theta_M)$ be such approximate solution corresponding to the initial data $(f_0,\Theta_0)$, which satisfy the hypothesis of Theorem \ref{global1}, then, in view of the Poincar\'e-like inequality in \eqref{lem:P}, we see that
$\mathscr{E}^0_0(f_M,\Theta_M)(t)$ is a Lyapunov functional, \emph{i.e.}, that
$$
\frac{d}{dt}\mathscr{E}^0_0(f_M,\Theta_M)(t)+\delta \mathscr{E}^0_0(f_M,\Theta_M)(0)\leq 0 \qquad \mbox{for all}\quad 0<t<\infty.
$$
The latter implies the exponential decay towards the equilibrium $(f^*,\Theta^*)=(0,0)$ for the global weak solution (see Subsection \ref{ss:decay}). Eventually, in Subsection \ref{ss:unique} it is shown that a global weak solution in 
$$
\left(L^1(0,T;\dot A^2(\TT))\right)^2
$$ is unique.

\medskip

\subsection{\emph{A priori} estimates in $\dot A^0(\TT)$:} 
\label{ss:e}
Let $T\in (0,\infty]$ and 
$$
(f,\Theta)\in \left( C^1([0,T), \dot A(\TT))\right)^2
$$ 
be a local solution of \eqref{eq:systemlinearized} with initial data 
$$
(f_0,\Theta_0)\in \left(\dot A(\TT)\right)^2
$$ 
satisfying the condition
\[
	\mathcal{E}_0^0(f_0,\Theta_0)<\min\{\bh,\bG\}.
\]
By continuity there exists $0<t_*\leq T$ such that the above estimate is satisfied  on $[0,t_*)$, that is
\begin{equation}\label{eq:t}
		\mathcal{E}_0^0(f,\Theta)(t)<\min\{\bh,\bG\}\qquad\mbox{for all}\quad t\in [0,t_*).
\end{equation}
Let us moreover assume that $t_*$ is the maximal time such that \eqref{eq:t} holds true on $[0,t_*)$. Hence, either
\begin{equation}\label{eq:O1}
\mathcal{E}_0^0(f,\Theta)(t_*)=\min\{\bh,\bG\}
\end{equation}
or $t_*=T$.
Notice that the restriction \eqref{restriction1} holds true on the time interval $[0,t_*)$. In the following let $t\in [0,t_*)$.
We compute
\begin{align*}
\pat|\hat{f}(k)| &=\frac{\text{Re}\left(\bar{\hat{f}}(k)\pat\hat{f}(k)\right)}{|\hat{f}(k)|},
\end{align*}
so, using that $\mathcal{S}=0$,
\begin{align}
\frac{d}{dt}\|f\|_{\dot A^0}&\leq\frac{\bh^2}{2}\|\Theta\|_{\dot{A}^2} +\left(\frac{\mathcal{A}}{\bh}-\frac{\mathcal{G}}{3}\bh^3\right)\| f\|_{\dot{A}^2}+\|N_1\|_{\dot A^0}+\|N_2\|_{\dot A^0}+\|N_4\|_{\dot A^0}\label{eq:1},\\
\frac{d}{dt}\|\Theta\|_{\dot A^0}&\leq-\left(\bh \bG+\mathcal{D}\right)\|\Theta\|_{\dot{A}^2}+\left|\frac{3\mathcal{A}\bG}{2\bh^2}-\bG\bh^2\frac{\mathcal{G}}{2}\right| \|f\|_{\dot{A}^2}+\|N_5\|_{\dot A^0}+\|N_6\|_{\dot A^0}+\|N_8\|_{\dot A^0}\label{eq:2}
\end{align}
We recall that
\begin{align*}
\widehat{uv}(k)=\sum_{j\in\ZZ}\hat{u}(j)\hat{v}(k-j)
\end{align*}
and by the hypothesis of the theorem
\[
	\|f\|_{\dot A^0}\leq \bh,\qquad \|\Theta\|_{\dot A^0}\leq \bG.
\]
Using the algebra property of the Wiener space $A^s(\TT)$, \emph{i.e.}
$$
\|fg\|_{\dot{A}^s}\leq 2^{s}\|f\|_{\dot{A}^s(\TT)}\|g\|_{\dot{A}^s(\TT)},\qquad \mbox{for all}\quad \, f,g\in \dot A^s(\TT),\quad s\in\NN,
$$
together with the interpolation inequality
$$
\|f\|_{\dot{A}^{s\theta}(\TT)}\leq \|f\|_{\dot A^0(\TT)}^{1-\theta}\|f\|_{\dot{A}^s(\TT)}^{\theta},\quad \mbox{for all} \quad 0<\theta<1,\quad s\geq0,
$$
the contribution of the nonlinear terms $N_k$, $k=1,2,$ and $4$, can be estimated as
\begin{align*}
\|N_1\|_{\dot A^0}&\leq \left(\frac{\|f\|_{\dot A^0}^2}{2}+\|f\|_{\dot A^0}\bh\right)\|\Theta\|_{\dot{A}^2}+\left(\|f\|_{\dot{A}^1}\|f\|_{\dot A^0}+\|f\|_{\dot{A}^1}\bh\right)\|\Theta\|_{\dot{A}^1}\\
&\leq \left(\frac{3}{2}\|f\|_{\dot A^0}\bh+\bh\|\Theta\|_{\dot A^0}\right)\|\Theta\|_{\dot{A}^2}+\bh\|f\|_{\dot A^0}\|f\|_{\dot{A}^2}\\
&\leq \mathscr{E}_0^0(f,\Theta) \frac{5}{2}\bh\|\Theta\|_{\dot{A}^2}+\bh\mathscr{E}_0^0(f,\Theta)\|f\|_{\dot{A}^2}\,,
\end{align*}
\begin{align*}
\|N_2\|_{\dot A^0}&\leq\frac{\mathcal{G}}{3}\bigg[\left(3\bh^2 \|f\|_{\dot A^0}+3\|f\|_{\dot A^0}^2\bh+\|f\|_{\dot A^0}^3\right)\|f\|_{\dot{A}^2}\\
&\quad+\left(3\bh^2 \|f\|_{A^1}+6\|f\|_{\dot{A}^1}\|f\|_{\dot A^0}\bh+3\|f\|_{\dot A^0}^2\|f\|_{\dot{A}^1}\right)\|f\|_{\dot{A}^1}\bigg]\\
&\leq\frac{19\bh^2\mathcal{ G}}{3} \|f\|_{\dot A^0}\|f\|_{\dot{A}^2}\\
&\leq \frac{19\bh^2\mathcal{ G}}{3} \mathscr{E}_0^0(f,\Theta)\|f\|_{\dot{A}^2}\,,
\end{align*}
and
\begin{align*}
\|N_4\|_{\dot A^0}&\leq\mathcal{A} \bigg{[}\left(\frac{\|f\|_{\dot A^0}}{\bh}\frac{\|f\|_{\dot{A}^2}}{\bh} +\left(\frac{\|f\|_{\dot{A}^1}}{\bh}\right)^2\right)\sum_{j=1}^\infty\left(\frac{\|f\|_{\dot A^0}}{\bh}\right)^{j-1}\\
&\quad +\frac{\|f\|_{\dot A^0}}{\bh}\left(\frac{\|f\|_{\dot{A}^1}}{\bh}\right)^2\sum_{j=1}^\infty j\left(\frac{\|f\|_{\dot A^0}}{\bh}\right)^{j-1} \bigg{]}\\
&\leq\mathcal{A} \bigg{[}\left(\frac{\|f\|_{\dot A^0}}{\bh}\frac{\|f\|_{\dot{A}^2}}{\bh} +\left(\frac{\|f\|_{\dot{A}^1}}{\bh}\right)^2\right)\left(1-\frac{\|f\|_{\dot A^0}}{\bh}\right)^{-1}\\
&\quad+\frac{\|f\|_{\dot A^0}}{\bh}\left(\frac{\|f\|_{\dot{A}^1}}{\bh}\right)^2\left(1-\frac{\|f\|_{\dot A^0}}{\bh}\right)^{-2} \bigg{]}\\
&\leq\mathcal{A} \frac{\|f\|_{\dot A^0}}{\bh}\frac{\|f\|_{\dot{A}^2}}{\bh}\left(1-\frac{\|f\|_{\dot A^0}}{\bh}\right)^{-1}\bigg{[}2+\left(1-\frac{\|f\|_{\dot A^0}}{\bh}\right)^{-1} \bigg{]}\\
&\leq \mathcal{A} \frac{\mathscr{E}_0^0(f,\Theta)}{\bh}\frac{\|f\|_{\dot{A}^2}}{\bh}\left(1-\frac{\|f\|_{\dot A^0}}{\bh}\right)^{-1}\bigg{[}2+\left(1-\frac{\|f\|_{\dot A^0}}{\bh}\right)^{-1} \bigg{]} ,
\end{align*}
where we used the convergence of the geometric series $\sum_{j=1}^\infty r^{j-1}$ for $|r|<1$, and
\[
	\sum_{j=1}^\infty jr^{j-1}=\partial_r \sum_{j=1}^{\infty}r^j = \frac{1}{(1-r)^2},\qquad |r|<1.
\]
Similarly,
the nonlinear terms $N_j$, $j=5,6,$ and $8$ are bounded by
\begin{align*}
\|N_5\|_{\dot A^0}
&\leq \left(\bG\|f\|_{\dot A^0}+\left(4\bh+\bG\right)\|\Theta\|_{\dot A^0} \right)\|\Theta\|_{\dot{A}^2}+ \bG\|f\|_{\dot A^0}\|f\|_{\dot{A}^2}\\
&\leq \mathscr{E}_0^0(f,\Theta)\left(4\bh+2\bG \right)\|\Theta\|_{\dot{A}^2}+ \bG\mathscr{E}_0^0(f,\Theta)\|f\|_{\dot{A}^2}\,,\end{align*}
\begin{align*}
\|N_6\|_{\dot A^0}&\leq\mathcal{G}\left[\left(4\bG\bh \|f\|_{\dot A^0}+4\|\Theta\|_{\dot A^0}\bh^2\right)\|f\|_{\dot{A}^2}+\bh^2\left(\|\Theta\|_{\dot A^0}\|\Theta\|_{\dot{A}^2}+\|f\|_{\dot A^0}\|f\|_{\dot{A}^2}\right)\right]\\
&\leq \mathcal{G}\mathscr{E}_0^0(f,\Theta)\left[\left(4\bG\bh +4\bh^2\right)\|f\|_{\dot{A}^2}+\bh^2\left(\|\Theta\|_{\dot{A}^2}+\|f\|_{\dot{A}^2}\right)\right]\,,
\end{align*}
\begin{align*}
\|N_8\|_{\dot A^0}&\leq \frac{3\mathcal{A}}{2\bh^4}\bigg{[}\left(3\|f\|_{\dot A^0}\bG\bh+\|\Theta\|_{\dot A^0} \bh^2\right)\|f\|_{\dot{A}^2}\\
&\quad+\left(4\|f\|_{\dot A^0}\bG\bh\|f\|_{\dot{A}^2}+\frac{\bh^2}{2}\left[\|\Theta\|_{\dot A^0}\|\Theta\|_{\dot{A}^2}+\|f\|_{\dot A^0}\|f\|_{\dot{A}^2}\right] \right)\bigg{]}\sum_{j=1}^\infty j\left(\frac{\|f\|_{\dot A^0}}{\bh}\right)^{j-1}\nonumber\\
&\quad+ \frac{3\mathcal{A}}{\bh^3}\left[2\bG\|f\|_{\dot A^0}\|f\|_{\dot{A}^2}\right]\sum_{j=2}^\infty j(j-1)\left(\frac{\|f\|_{\dot A^0}}{\bh}\right)^{j-2}\\
&\leq \frac{3\mathcal{A}}{2\bh^3}\bigg{[}\left(3\|f\|_{\dot A^0}\bG+\|\Theta\|_{\dot A^0} \bh\right)\|f\|_{\dot{A}^2}\\
&\quad+\left(4\|f\|_{\dot A^0}\bG\|f\|_{\dot{A}^2}+\frac{\bh}{2}\left[\|\Theta\|_{\dot A^0}\|\Theta\|_{\dot{A}^2}+\|f\|_{\dot A^0}\|f\|_{\dot{A}^2}\right] \right)\bigg{]}\left(1-\frac{\|f\|_{\dot A^0}}{\bh}\right)^{-2}\nonumber\\
&\quad+ \frac{3\mathcal{A}}{\bh^3}\left[2\bG\|f\|_{\dot A^0}\|f\|_{\dot{A}^2}\right]\left(1-\frac{\|f\|_{\dot A^0}}{\bh}\right)^{-3}.
\end{align*}
Grouping terms, we find that
\begin{align*}
\|N_8\|_{\dot A^0}&\leq \frac{3\mathcal{A}}{2\bh^3}\left(1-\frac{\|f\|_{\dot A^0}}{\bh}\right)^{-2}\Big\{7\bG\|f\|_{\dot A^0}\|f\|_{\dot A^2} + \frac{1}{2}\bh[\|f\|_{\dot A^0}\|f\|_{\dot A^2}+{\frac{\bh}{2}}\|\Theta\|_{\dot A^0}\|\Theta\|_{\dot A^2}] \\
  &\qquad + \bh \|\Theta\|_{\dot A^0}\|f\|_{\dot A^2}+ \left(1-\frac{\|f\|_{\dot A^0}}{\bh}\right)^{-1}4\bG\|f\|_{\dot A^0}\|f\|_{\dot A^2}\Big\}\\
&\leq \frac{3\mathcal{A}}{2\bh^3}\mathscr{E}_0^0(f,\Theta)\left(1-\frac{\|f\|_{\dot A^0}}{\bh}\right)^{-2}\Big\{7\bG\|f\|_{\dot A^2} + \frac{1}{2}\bh[\|f\|_{\dot A^2}+{\frac{\bh}{2}}\|\Theta\|_{\dot A^2}] \\
  &\qquad + \bh \|f\|_{\dot A^2}+ \left(1-\frac{\|f\|_{\dot A^0}}{\bh}\right)^{-1}4\bG\|f\|_{\dot A^2}\Big\}.
\end{align*}

We recall the definition \eqref{constant1} and \eqref{constant2}, then we add equations \eqref{eq:1} and \eqref{eq:2}, use the previous estimates for $N_j$, $j=1,2,4,5,6,$ and 8 and  obtain that
\begin{align}\label{eq:energyE}
\frac{d}{dt}\mathscr{E}_0^0(f,\Theta)(t)&\leq -\left(\mathfrak{C}_1-\Lambda_1(t)\mathscr{E}_0^0(f,\Theta)(t)\right)\|f(t)\|_{\dot{A}^2}-\left(\mathfrak{C}_2-\Lambda_2(t)\mathscr{E}_0^0(f,\Theta)(t)\right)\|\Theta(t)\|_{\dot{A}^2}
\end{align}
 for $t\in[0,t_*)$,
with $\Lambda_1$ and $\Lambda_2$ defined in \eqref{eq:Lambda_1} and \eqref{eq:Lambda_2},  respectively. 
Using the hypothesis of the Theorem \ref{global1}, we have that
\begin{align*}
\frac{d}{dt}\mathscr{E}^0_0(f,\Theta)(t)\bigg{|}_{t=0}&<0,
\end{align*}
so, there exists a time $0<t_0\leq t_*$ such that
\begin{equation}\label{eq:estimate}
\mathscr{E}_0^0(f,\Theta)(t)\leq \mathscr{E}^0_0(f_0,\Theta_0)\qquad \mbox{for all}\quad t\in [0,t_0].
\end{equation}
Let us assume that $t_0$ is the maximal times such that \eqref{eq:estimate} holds true on $[0,t_0]$
We want to propagate this decay for all times, that is, we aim to show that $t_0=t_*$, which in turn implies that $t_0=t_*=T$, by \eqref{eq:O1}. Let us emphasize that $\mathscr{E}_0^0(f,\Theta)(t)\leq \mathscr{E}^0_0(f_0,\Theta_0)$ for $t\in [0,t_0]$ guarantees that
\begin{equation}\label{eq:Ldecay}
\Lambda_j(t)\leq \Lambda_j(0),\qquad j=1,2 \qquad\mbox{for all}\quad t\in [0,t_0].
\end{equation}
Thereby, for any $t\in [0,t_0]$ we have that
$$
\mathfrak{C}_1-\Lambda_1(t)\mathscr{E}^0_0(f,\Theta)(t)\geq \mathfrak{C}_1-\Lambda_1(0)\mathscr{E}^0_0(f_0,\Theta_0)=\gamma_1>0,
$$
$$
\mathfrak{C}_2-\Lambda_2(t)\mathscr{E}^0_0(f,\Theta)(t)\geq \mathfrak{C}_2-\Lambda_2(0)\mathscr{E}^0_0(f_0,\Theta_0)=\gamma_2>0.
$$
In particular, we obtain that 
\begin{align*}
\frac{d}{dt}\mathscr{E}_0^0(f,\Theta)(t)&< 0\qquad\mbox{for all}\quad t\in [0,t_0).
\end{align*}
Assume that $t_0<t_*$. By continuity, we deduce that $\mathscr{E}_0^0(f,\Theta)(t_0)=\mathscr{E}_0^0(f,\Theta)(0)$, but this implies that, again, 
\begin{align*}
\frac{d}{dt}\mathscr{E}_0^0(f,\Theta)(t)\bigg{|}_{t=t_0}&\leq 0,
\end{align*}
and that contradicts the assumption $t_0<t_*$. Thus, we have shown that in fact
\begin{equation}\label{eq:Edecay}
\mathscr{E}_0^0(f,\Theta)(t)\leq \mathscr{E}_0^0(f_0,\Theta_0)\qquad \mbox{for all}\quad t\in[0,t_*)
\end{equation}
and thereby $t_*=T$,
in view of \eqref{eq:O1}.
Then,
\begin{equation}\label{eq:E}
\frac{d}{dt}\mathscr{E}_0^0(f,\Theta)(t)\leq -\gamma \mathscr{E}^2_2(f,\Theta)(t)\qquad \mbox{for all}\quad t\in[0,T),
\end{equation}
where $\delta:=\min\{\gamma_1,\gamma_2\}$. Eventually, the energy estimate for the gravity driven equation ($\mathcal S =0$) reads
\[
	\mathscr{E}^0_0(f,\Theta)(t)+\delta\int_0^t \mathscr{E}^2_2(f,\Theta)(\tau)\, d\tau \leq \mathscr{E}^0_0(f_0,\Theta_0)
\]
and
\[
	\int_0^t \|\partial_t f(\tau)\|_{\dot A^0}+\|\partial_t \Theta(\tau)\|_{\dot A^0}\, d\tau\leq c,
\]
for all $t\in[0,T)$, where $c>0$ is a constant depending on the initial data.

\medskip

\subsection{Existence of global weak solutions}\label{ss:gl} 
We use a standard Galerkin approximation to obtain in the limit a global weak solution of \eqref{eq:systemlinearized} where surface tension effects are neglected, \emph{i.~e.} $\mathcal{S}=0$. Let us fix $M\in \ZZ^+$. Set
\[
	f_M(t,x):=\sum_{|k|\leq M} \hat f(t,k)e^{ikx}\qquad \mbox{and}\qquad 	\Theta_M(t,x):=\sum_{|k|\leq M} \hat \theta(t,k)e^{ikx}
\]
and the initial data
\[
	f_M(0,x):=\sum_{|k|\leq M} \hat f_0(k)e^{ikx}\qquad \mbox{and}\qquad 	\Theta_M(0,x):=\sum_{|k|\leq M} \hat \theta_0(k)e^{ikx}
\]
to coincide with the Fourier truncation of the $f_0$ and $\Theta_0$, respectively. Recall that the convergence of the Fourier series of the initial data is guaranteed by the assumption that $f_0,\Theta_0\in \dot A^0(\TT)$. We consider the Galerkin approximated problems:
	\begin{alignat*}{2}
	\partial_t f_M-\frac{\bh^2}{2}\partial_x^2 \Theta_M +\left(\frac{\mathcal{A}}{\bh}-\frac{\mathcal{G}}{3}\bh^3\right)\partial_x^2 f_M+\frac{\mathcal{S}}{3}\bh^3\partial_x^4 f_M&=\sum_{j=1,2,4} N_j^M,\; && \text{ in }(0,T)\times\TT\\
	\partial_t \Theta_M -\left(\bh\bG+\mathcal{D}\right)\pax^2\Theta_M+\left(\frac{3\mathcal{A}\bG}{2\bh^2}-\frac{\mathcal{G}}{2}\bG\bh^2\right) \pax^2f_M+\frac{\mathcal{S}}{2}\bG\bh^2 \pax^4f_M &= \sum_{j=5,6,8} N_j^M,\;&& \text{ in }(0,T)\times\TT
	\end{alignat*}
where the nonlinearities $N_j^M$ are given by

\begin{align*}
N_1^M&=\pax P_M\left[\left(\frac{f^2}{2}+f\bh\right)\partial_x \Theta\right],\\
N_2^M&=\partial_x P_M \left[\frac{\mathcal{G}}{3}\left(3\bh^2 f+3f^2\bh+f^3\right)\partial_x f\right],\\
N_4^M&=\mathcal{A} P_M \bigg{[}\left(\frac{f}{\bh^2}\partial_x^2 f +\left(\frac{\partial_x f}{\bh}\right)^2\right)\sum_{j=1}^M (-1)^{j+1}\left(\frac{f}{\bh}\right)^{j-1}-\frac{f}{\bh^3}(\partial_x f)^2\sum_{j=1}^M j(-1)^{j+1}\left(\frac{f}{\bh}\right)^{j-1} \bigg{]},\\
N_5^M&=\pax P_M\left[\left(\bG f+\Theta\bh+\Theta f\right)\partial_x \Theta\right],\\
N_6^M&=\pax P_M\left[\frac{\mathcal{G}}{2}\left(\bG f^2+2\bG\bh f+\Theta\bh^2+\Theta f^2+2\Theta\bh f\right)\partial_x f\right],\\
N_8^M&=\frac{3\mathcal{A}}{2\bh^4} P_M\left[\left(2f\bG\bh+f^2\bG-\Theta \bh^2\right)\partial_x^2 f+\left(2\pax f\bG\bh+2f\pax f\bG-\pax \Theta \bh^2\right)\partial_x f\right]\sum_{j=1}^M j(-1)^{j+1}\left(\frac{f}{\bh}\right)^{j-1}\nonumber\\
&\quad- \frac{3\mathcal{A}}{\bh^4} P_M\left[\left(2f\bG\bh+f^2\bG-\Theta \bh^2\right)\frac{(\partial_x f)^2}{\bh}\right]\frac{1}{2}\sum_{j=2}^M j(j-1)(-1)^{j}\left(\frac{f}{\bh}\right)^{j-2}.
\end{align*}
Here, the operator $P_M$ denotes the Fourier truncation operator
\[
	P_M g(x)=\sum_{|k|\leq M} \hat g(k)e^{ikx}\qquad\mbox{for any}\quad g\in \dot A^0(\TT).
\]
The Picard--Lindel\"of Theorem ensures the existence of classical solutions 
\[
f_M,\Theta_M\in C^1([0,T_M);C^\infty(\TT)),
\]
 where $T_M>0$ is the maximal existence time.  Furthermore, the approximated problems provide the same \emph{a priori} bounds as in the previous sections.  Consequently the Galerkin solutions  $(f_M,\Theta_M)_{M\in \NN}$ exist globally and for any $T>0$ we have the bounds:
\begin{equation}\label{boundedness}
	(f_M,\Theta_M)_{M\in\NN} \quad \mbox{is uniformly bounded in}\quad \left( L^\infty\left(0,T;\dot A^0(\TT)\right)\cap L^1\left(0,T; \dot A^2(\TT) \right)\right)^2
\end{equation}
and
\begin{equation}\label{boundednesst}
(\partial_t f_M,\partial_t\Theta_M)_{M\in\NN} \quad \mbox{is uniformly bounded in}\quad \left( L^1\left(0,T; \dot A^0(\TT) \right)\right)^2.
\end{equation}
Following the lines in \cite{BG18}, the above uniform regularities of the Galerkin approximations $(f_M,\Theta_M)_{M\in\NN} $ guarantee the existence of a weakly convergent subsequences (not relabeled) such that 
\begin{align*}
(f_M,\Theta_M)\overset{*}{\rightharpoonup}  (f,\Theta)\text {  in $\left(L^{\infty}(0,T;L^\infty(\TT))\right)^2$}.
\end{align*}
Similarly, using interpolation in Wiener spaces, the finite measure of the spatial domain and the fact that $\dot A^r(\TT)\subset \dot W^{r,\infty}(\TT)$ for any $r\geq 0$, we obtain the existence of a subsequence (not relabeled), such that
\begin{align}\label{eq:Linf1}
(f_M,\Theta_M)\rightharpoonup  (f,\Theta)\text {   in $\left(L^{\frac{2}{r}}(0,T;\dot W^{r,p}(\TT))\right)^2$},\qquad 0\leq r<2,\;\;1\leq p<\infty.
\end{align}
From the previous fact we can infer that actually
$$
(f,\Theta) \in \left(L^{\frac{2}{r}}(0,T;\dot W^{r,\infty}(\TT))\right)^2,\qquad 0\leq r<2.
$$
Furthermore \eqref{boundedness}, \eqref{boundednesst} imply that
$$
	(f_M,\Theta_M)_{M\in\NN} \text{ is uniformly bounded in }   \left(L^2(0,T;H^{1}(\TT))\right)^2,
$$
which is due to $\dot A^s(\TT)\subset \dot H^s(\TT)$ for any $s \geq 0$ and an interpolation inequality for fractional Sobolev spaces. Eventually, as a consequence of \eqref{boundedness},\eqref{boundednesst} and a  compactness argument as in \cite[Corollary 4]{simon1986compact}, we obtain (up to a subsequence) that
$$
(f_M,\Theta_M)\rightarrow (f,\Theta)\qquad \text{ in $\left(L^1(0,T;{C}^{1+s}(\TT))\right)^2$}, \quad \,0\leq s<\frac{1}{2}.
$$
Passing to the limit in the weak formulation of the Galerkin approximation yields the existence of a global weak solution of \eqref{eq:systemlinearized} in the sense of Definition \ref{defi1}.

\subsection{Exponential trend to equilibrium}\label{ss:decay} Using \eqref{eq:energyE} and the definition of $\delta>0$ in \eqref{delta}, we have that
\begin{align}\label{eq:expdecay}
\frac{d}{dt}\mathscr{E}^0_0(f_M,\Theta_M)(t)\leq -\delta\mathscr{E}^2_2(f_M,\Theta_M)(t) \qquad \mbox{for all}\quad t\geq0.
\end{align}
Using the Poincar\'e-like inequality \eqref{lem:P}, we also conclude that 
$$
\frac{d}{dt}\mathscr{E}^0_0(f_M,\Theta_M)(t)\leq -\delta\mathscr{E}^0_0(f_M,\Theta_M)(t) \qquad \mbox{for all}\quad t\geq0,
$$ 
 which in turn implies the exponential decay towards the equilibrium:
\[
\mathscr{E}^0_0(f_M,\Theta_M)(t)\leq \mathscr{E}_0^0(f_0,\Theta_0)e^{-\delta t}\qquad \mbox{for all}\quad t\geq0.
\]
Using \eqref{eq:Linf1} and the lower semi-continuity of the weak$-*$ convergence, we have that
\begin{equation*}
\|f(t)\|_{L^\infty(\TT)}+\|\Theta(t)\|_{L^\infty(\TT)}\leq \liminf_{k\to \infty}\mathscr{E}_0^0(f_M,\Theta_M)(t)\leq 	\mathscr{E}_0^0(f_0,\Theta_0)e^{-\delta t}.
\end{equation*}

\subsection{Uniqueness}
\label{ss:unique}
The proof follows a standard contradiction argument. For the sake of brevity, we only sketch the idea. Assume that there exist two different solutions $(f_1,\Theta_1)$ and $(f_2,\Theta_2)$ starting from the same initial data 
$$
(f_0,\Theta_0)\in \left(\dot A^0(\TT)\right)^2.
$$ 
Assume also that these solutions satisfy 
$$
(f_i,\Theta_i)\in \left( L^{1}(0,T;\dot A^2(\TT))\right)^2.
$$ Using the smallness of the initial data, the same estimates as in Subsection \ref{ss:e} yield that
\begin{multline*}
\frac{d}{dt}\left(\|f_1-f_2\|_{\dot A^0}+\|\Theta_1-\Theta_2\|_{\dot A^0}\right)\\
\leq C\left[\|f_1-f_2\|_{\dot A^0}+\|\Theta_1-\Theta_2\|_{\dot A^0}\right]\left(\|f_1\|_{\dot A^2}+\|f_2\|_{\dot A^2}+\|\Theta_1\|_{\dot A^2}+\|\Theta_2\|_{\dot A^2}+1\right).
\end{multline*} 

Now the statement is a consequence of Gronwall's  inequality and the fact that $(f_1,\Theta_1)(0)=(f_2,\Theta_2)(0)$.

\section{Proof of Theorem \ref{global2}: Global existence and decay when $\mathcal{S}>0$}
\label{S:T2}
The proof essentially follows the arguments in the previous section, the main difference relying in the fact that for $\mathcal S>0$ we have that \eqref{eq:systemlinearized}  is a system of mixed orders. Thereby the energy estimates require some additional investigation. The existence of local solutions of the approximated Galerkin systems are straightforward due to Picard--Lindel\"of's theorem (see Subsection \ref{ss:gl}). The energy estimates then ensure that the approximated solutions exist globally. Furthermore, the energy estimates provide the necessary a priori bounds to pass to the limit in the Galerkin approximation; thereby guaranteeing the existence of a global weak solution in the sense of Definition \ref{defi1}. In view of the weak lower semicontinuity of the norm, the global weak solution inherits the energy estimates for the approximated solutions and we can conclude the regularity and exponential decay of the solution.

\subsection{\emph{A priori} estimates in $\dot A^0(\TT)$:} 
To perform our energy estimates let us assume that there exists a local solution 
$$
(f,\Theta)\in \left(C^1([0,T), \dot A^0(\TT)) \right)^2
$$ 
of the fourth-order system \eqref{eq:systemlinearized} where $\mathcal{S}>0$, corresponding to initial data $(f_0,\Theta_0)$, satisfying 
\[
	\mathscr{E}_0^0(f_0,\Theta_0)<\min\{\bh,\bG\}.
\]

Similar as before, we have that
\begin{align}
\frac{d}{dt}\|f\|_{\dot A^0}&\leq\frac{\bh^2}{2}\|\Theta\|_{\dot{A}^2} +\left(\frac{\mathcal{A}}{\bh}-\frac{\mathcal{G}}{3}\bh^3\right)\| f\|_{\dot{A}^2} - \frac{\mathcal{S}}{3}\bh^3\|f\|_{\dot A^4}+\sum_{j=1}^4\|N_j\|_{\dot A^0}\label{eq:1b},\\
\frac{d}{dt}\|\Theta\|_{\dot A^0}&\leq-\left(\bh \bG+\mathcal{D}\right)\|\Theta\|_{\dot{A}^2}+\left|\frac{3\mathcal{A}\bG}{2\bh^2}-\bG\bh^2\frac{\mathcal{G}}{2}\right|\|f\|_{\dot A^2}+\frac{\mathcal{S}}{2}\bh^2\bG\|f\|_{\dot A^4} +\sum_{j=5}^8\|N_j\|_{\dot A^0}.\label{eq:2b}
\end{align}
Keeping in mind the definitions of $\mathfrak{C_i}$, $i=1,2,3,$ in \eqref{constant1}, \eqref{constant2}, and \eqref{constant3}, we take the sum of the two inequalities above and obtain that

$$
\frac{d}{dt}\mathscr{E}_0^0(f,\Theta)+\mathfrak{C}_1\|f\|_{\dot{A}^2}+\mathfrak{C}_2\|\Theta\|_{\dot{A}^2}+\mathfrak{C}_3\|f\|_{\dot{A}^4}\leq \sum_{j=1}^8\|N_j\|_{\dot A^0}.
$$
Recalling \eqref{eq:energyE}, we have that
\begin{align*}
\frac{d}{dt}&\mathscr{E}^0_0(f,\Theta)+\mathfrak{C}_3\|f\|_{\dot{A}^4}\\
&\leq \|N_3\|_{\dot A^0}+\|N_7\|_{\dot A^0} -\left(\mathfrak{C}_1-\Lambda_1(t)\mathscr{E}^0_0(t)\right)\|f(t)\|_{\dot{A}^2}-\left(\mathfrak{C}_2-\Lambda_2(t)\mathscr{E}^0_0(t)\right)\|\Theta(t)\|_{\dot{A}^2},
\end{align*}
where $\Lambda_1$ and $\Lambda_2$ are defined in \eqref{eq:Lambda_1} and \eqref{eq:Lambda_2}, respectively.
Thus we are left to estimate the remaining terms $\|N_3\|_{\dot A^0}$ and $\|N_7\|_{\dot A^0}$. Similarly as before, using the interpolation inequality in Wiener spaces, we estimate
\begin{align*}
\|N_3\|_{\dot A^0}&\leq \frac{19}{3}\mathcal{S}\bh^2\mathscr{E}_0^0(f,\Theta)\|f\|_{\dot{A}^4},
\end{align*}
and
\begin{align*}
\|N_7\|_{\dot A^0}&\leq \frac{\mathcal{S}}{2}\left(6\bG\bh \|f\|_{\dot A^0}+\bh^2\|\Theta\|_{\dot A^0}\right)\|f\|_{\dot{A}^4}\\
&\quad +\frac{\mathcal{S}}{2}\left(8\bG \bh \|f\|_{A^1}+4\|\Theta\|_{A^1}\bh^2\right)\|f\|_{\dot{A}^3}.
\end{align*}
Note that by Young's inequality, we have that
$$
\|\Theta\|_{\dot{A}^1}\|f\|_{\dot{A}^3}\leq \|\Theta\|_{\dot{A}^1}\|f\|_{\dot A^0}^{\frac{1}{4}}\|f\|_{\dot{A}^4}^{\frac{3}{4}}\leq\|\Theta\|_{\dot{A}^1}\left(\frac{1}{4\varepsilon}\|f\|_{\dot A^0}+\frac{3\varepsilon}{4}\|f\|_{\dot{A}^4}\right), 
$$
for any $\varepsilon>0$. In particular, as $\Theta$ and $f$ are nonzero (provided the initial data are nontrivial), we can take
$$
\varepsilon=\frac{\mathscr{E}_0^0(f,\Theta)}{\|\Theta\|_{\dot{A}^1}}.
$$
Thus,
$$
\|\Theta\|_{\dot{A}^1}\|f\|_{\dot{A}^3}\leq\frac{\|\Theta\|_{\dot{A}^1}^2}{4\mathscr{E}_0^0(f,\Theta)}\|f\|_{\dot A^0}+\frac{3\mathscr{E}_0^0}{4}\|f\|_{\dot{A}^4}\leq \frac{\|\Theta\|_{\dot{A}^2}\mathscr{E}_0^0(f,\Theta)}{4}+\frac{3\mathscr{E}_0^0(f,\Theta)}{4}\|f\|_{\dot{A}^4}. 
$$
Then, we conclude that
	\[
		\|N_7\|_{\dot A^0}\leq \frac{\mathcal{S}}{2}\left( (14\bh\bG +4\bh^2)\mathscr{E}_0^0(f,\Theta) \|f\|_{\dot A^4}+\bh^2\mathscr{E}_0^0(f,\Theta) \|\Theta\|_{\dot A^2}\right).
	\]
Using \eqref{eq:Lambda_3}, we can group terms as follows:
\begin{align*}
\frac{d}{dt}\mathscr{E}^0_0(f,\Theta)(t)&\leq -\left(\mathfrak{C}_1-\Lambda_1(t)\mathscr{E}^0_0(f,\Theta)(t)\right)\|f(t)\|_{\dot{A}^2}-\left(\mathfrak{C}_2-\left(\Lambda_2(t)+\frac{\mathcal{S}\bh^2}{2}\right)\mathscr{E}^0_0(f,\Theta)(t)\right)\|\Theta(t)\|_{\dot{A}^2}\\
&\quad-\left(\mathfrak{C}_3-\Lambda_3\mathscr{E}^0_0(f,\Theta)(t)\right)\|f(t)\|_{\dot{A}^4}.
\end{align*}
Repeating the argument from Section \ref{S:T1}, we conclude  that
\begin{equation}\label{eq:Es}
\mathscr{E}^0_0(f,\Theta)(t)\leq \mathscr{E}^0_0(f_0,\Theta_0)e^{-\delta t} 
\end{equation}
where $\delta$ is defined in \eqref{delta}. The energy estimate for the capillary driven thin film ($\mathcal S >0$) reads
\begin{equation*}
\mathscr{E}^0_0(f,\Theta)(t)+ \int_0^t\delta \mathscr{E}^4_2(f,\Theta)(\tau)\, d\tau \leq \mathscr{E}^0_0(f_0,\Theta_0)
\end{equation*}
and
\[
\int_0^t \|\partial_t f(\tau)\|_{\dot A^0}+\|\partial_t \Theta(\tau)\|_{\dot A^0}\,d\tau \leq c,
\]
for all $t\in[0,T)$, where $c>0$ is a constant depending on the initial data.

\subsection{Existence of global weak solutions} In a similar way as in Subsection \ref{ss:gl}, we obtain the existence of a sequence of global Galerkin approximations $(f_M,\Theta_M)_{M\in\NN}$. Then, the energy estimates from above guarantee that for any $T>0$:
$$
(f_M)_{M\in\NN}\quad\mbox{is uniformly bounded in}\quad  L^\infty(0,T; A^0(\TT))\cap  L^1(0,T;A^4(\TT))\subset  L^2(0,T;H^2(\TT))
$$
and
$$
(\Theta_M)_{M\in\NN}\quad\mbox{is uniformly bounded in}\quad L^\infty(0,T; A^0(\TT))\cap L^1(0,T;A^2(\TT))\subset   L^2(0,T;H^1(\TT)).
$$
Moreover, the time derivatives satisfy
$$
(\partial_t f_M,\partial_t \Theta_M)_{M\in\NN}\quad\mbox{is uniformly bounded in}\quad \left(L^1(0,T;A^0(\TT))\right)^2
$$
for any $T>0$. Consequently, we obtain the existence of a subsequence (not relabeled) such that
\begin{align}\label{eq:Linf}
\begin{split}
&f_M\rightharpoonup  f\text {   in $L^{\frac{4}{s}}(0,T; W^{s,p}(\TT))$},\quad 0\leq s<4,\;\;1\leq p<\infty\\
&\Theta_M\rightharpoonup  \Theta\text {   in $L^{\frac{2}{r}}(0,T; W^{r,p}(\TT))$} ,\quad 0\leq r<2,\;\;1\leq p<\infty.
\end{split}
\end{align}
Moreover,
\begin{align*}
f_M&\rightarrow f\qquad \text{ in $L^1(0,T;C^{3+\alpha}(\TT))$}, \quad \,0\leq \alpha<\frac{1}{2},\\
f_M&\rightharpoonup f\qquad \text{ in $L^2(0,T, H^2(\TT))$}
\end{align*}
and 
\begin{align*}
\Theta_M&\rightarrow \Theta\qquad \text{ in $L^1(0,T;C^{1+\alpha}(\TT))$},\quad  \,0\leq \alpha<\frac{1}{2},\\
\Theta_M&\rightharpoonup \Theta\qquad \text{ in $L^2(0,T, H^1(\TT))$}.
\end{align*}
Equipped with these convergences we can pass to the limit in the weak formulation and conclude the global existence of a weak solution of \eqref{eq:systemlinearized} for $\mathcal S>0$ in the sense of Definition \ref{defi1}.

\subsection{Exponential trend to equilibrium} The proof follows the same ideas as in Subsection~\ref{ss:decay}.
The sequence of Galerkin approximation $(f_M, \Theta_M)_{M\in\NN}$ satisfies the energy estimates. In particular, we have the exponential decay towards the equilibrum in \eqref{eq:Es}:
 \[
 \mathscr{E}^0_0(f_M,\Theta_M)(t)\leq \mathscr{E}_0^0(f_0,\Theta_0)e^{-\delta t}\qquad \mbox{for all}\quad t\geq0,
 \]
 where $\delta>0$ is defined in \eqref{delta}.
 Using \eqref{eq:Linf} and the lower semi-continuity of the weak$-*$ convergence, we have that
 \begin{equation*}
 \|f(t)\|_{L^\infty(\TT)}+\|\Theta(t)\|_{L^\infty(\TT)}\leq \liminf_{k\to \infty}\mathscr{E}^0_0(f_M,\Theta_M)(t)\leq 	\mathscr{E}_0^0(f_0,\Theta_0)e^{-\delta t}
 \end{equation*}
for all $t\geq 0$, which proves the claim.

\subsection{Uniqueness} Also the proof for the conditional uniqueness is similar to the one in Subsection \ref{ss:unique}.
 Assume that there exist two different solutions $(f_1,\Theta_1)$ and $(f_2,\Theta_2)$ starting from the same initial data 
 $$
 (f_0,\Theta_0)\in \left(\dot A(\TT)\right)^2.
 $$ Moreover, we suppose that the solutions satisfy the additional regularity 
 $$
 (f_i,\Theta_i)\in \left( L^{1}(0,T;\dot A^4(\TT))\times L^1(0,T;\dot A^2(\TT))\right)^2
 $$ for $i=1,2$. Similar as for the energy estimate we compute that
 \begin{multline*}
 \frac{d}{dt}\left(\|f_1-f_2\|_{\dot A^0}+\|\Theta_1-\Theta_2\|_{\dot A^0}\right)\\
 \leq C \left[\|f_1-f_2\|_{\dot A^0}+\|\Theta_1-\Theta_2\|_{\dot A^0}\right]\\
 \times\left(\|f_1\|_{\dot A^2}+\|f_2\|_{\dot A^2}+\|f_1\|_{\dot A^4}+\|f_2\|_{\dot A^4}+\|\Theta_1\|_{\dot A^2}+\|\Theta_2\|_{\dot A^2}\right).
 \end{multline*} 
 Now the assertion follows by applying Gronwall's inequality and recalling that the solutions $(f_i,\Theta_i), i=1,2$ share the same initial data.
 
 \bigskip
 
\section*{Acknowledgements}The author R.G.B was partially supported by the LABEX MILYON (ANR-10-LABX-0070) of Universit\'e de Lyon, within the program ``Investissements d'Avenir'' (ANR-11-IDEX-0007) operated by the French National Research Agency (ANR).  The author G.B. gratefully acknowledges the funding by the Deutsche Forschungsgemeinschaft (DFG, German Research Foundation)- Project ID 258734477 - SFB 1173.
Part of this research was carried out while G.B. was supported by grant no. 250070 from the Research Council of Norway and during a short  stay of G.B. at Institut Camille Jordan under project DYFICOLTI, ANR-13-BS01-0003-01 support. 

\bigskip

\bibliographystyle{plain}

\begin{thebibliography}{10}

\bibitem{BGN03}
J.~W. Barrett, H.~Garcke, and R.~N{\"{u}}rnberg.
\newblock {Finite element approximation of surfactant spreading on a thin
  film}.
\newblock {\em SIAM J. Numer. Anal.}, 41(4):1427--1464, 2003.

\bibitem{BN04}
J.~W. Barrett and R.~N{\"{u}}rnberg.
\newblock {Convergence of a finite-element approximation of surfactant
  spreading on a thin film in the presence of van der Waals forces}.
\newblock {\em IMA J. Numer. Anal.}, 24(2):323--363, 2004.

\bibitem{Beretta1995}
E.~Beretta, M.~Bertsch, and R.~Dal~Passo.
\newblock {Nonnegative solutions of a fourth-order nonlinear degenerate
  parabolic equation}.
\newblock {\em Archive for Rational Mechanics and Analysis}, 129(2):175--200,
  1995.

\bibitem{Bernis1990}
F.~Bernis and A.~Friedman.
\newblock {Higher order nonlinear degenerate parabolic equations}.
\newblock {\em Journal of Differential Equations}, 83(1):179--206, 1990.

\bibitem{Bertozzi1996}
A.~L. Bertozzi and M.~C. Pugh.
\newblock The lubrication approximation for thin viscous films : regularity
  and long-time behavior of weak solutions.
\newblock {\em Comm. Pure Appl. Math.}, 49(2):85--123, 1996.

\bibitem{BG}
M.~S. Borgas and J.~B. Grotberg.
\newblock {Monolayer flow on a thin film}.
\newblock {\em J. Fluid Mech.}, 193:151--170, 1988.

\bibitem{B1}
G.~Bruell.
\newblock {Modeling and analysis of a two-phase thin film model with insoluble
  surfactant}.
\newblock {\em Nonlinear Anal. Real World Appl.}, 27:124--145, 2016.

\bibitem{B2}
G.~Bruell.
\newblock {Weak solutions to a two-phase thin film model with insoluble
  surfactant driven by capillary effects}.
\newblock {\em Journal of Evolution Equations}, 17(4):1341--1379, 2017.

\bibitem{bruellwell}
G. Bruell.
\newblock Well-posedness and stability for a mixed order system arising in thin
  film equations with surfactant.
\newblock {\em to appear in Mathematische Nachrichten (2019)}.

\bibitem{BG18}
G. Bruell and R. Granero-Belinch{\'o}n.
\newblock On the thin film Muskat and the thin film Stokes equations.
\newblock {\em Journal of Mathematical Fluid Mechanics}, 21(2):33, 2019.

\bibitem{burczak2016generalized}
J. Burczak and R. Granero-Belinch{\'o}n.
\newblock On a generalized doubly parabolic Keller--Segel system in one spatial
  dimension.
\newblock {\em Mathematical Models and Methods in Applied Sciences},
  26(01):111--160, 2016.

\bibitem{CT13}
M.~Chugunova and R.~M. Taranets.
\newblock {Nonnegative weak solutions for a degenerate system modeling the
  spreading of surfactant on thin films}.
\newblock {\em Appl. Math. Res. Express. AMRX}, (1):102--126, 2013.

\bibitem{constantin1993droplet}
P. Constantin, T. Dupont, R. Goldstein, L. Kadanoff, M. Shelley, and S. Zhou.
\newblock Droplet breakup in a model of the Hele-Shaw cell.
\newblock {\em Physical Review E}, 47(6):4169, 1993.

\bibitem{gan2}
P. Constantin, D. Córdoba, F. Gancedo, L. Rodríguez-Piazza and R. Strain.
\newblock {On the Muskat problem: global in time results in 2D and 3D}.
\newblock {\em American Journal of Mathematics}, 138(6):1455--1494, 2016.

\bibitem{gan3}
P. Constantin, D. Córdoba, F. Gancedo and R. Strain.
\newblock {On the global existence for the Muskat problem}.
\newblock {\em J. Eur. Math. Soc.}, 15(1), 201--227, 2013.

\bibitem{ECM}
B.~D. Edmonstone, R.~V. Craster, and O.~K. Matar.
\newblock {Surfactant-induced fingering phenomena beyond the critical micelle
  concentration}.
\newblock {\em J. Fluid Mech.}, 564:105--138, 2006.

\bibitem{EHLW11}
J.~Escher, M.~Hillairet, Ph. Lauren{\c{c}}ot, and Ch. Walker.
\newblock {Global weak solutions for a degenerate parabolic system modeling the
  spreading of insoluble surfactant}.
\newblock {\em Indiana Univ. Math. J.}, 60(6):1975--2019, 2011.

\bibitem{EHLW1}
J.~Escher, M.~Hillairet, Ph. Lauren{\c{c}}ot, and Ch. Walker.
\newblock {Thin film equations with soluble surfactant and gravity: modeling
  and stability of steady states}.
\newblock {\em Math. Nachr.}, 285(2-3):210--222, 2012.

\bibitem{EHLW4th}
J.~Escher, M.~Hillairet, Ph. Lauren{\c{c}}ot, and Ch. Walker.
\newblock {Weak solutions to a thin film model with capillary effects and
  insoluble surfactant}.
\newblock {\em Nonlinearity}, 25(9):2423--2441, 2012.

\bibitem{EHLWtrav}
J.~Escher, M.~Hillairet, Ph. Lauren{\c{c}}ot, and Ch. Walker.
\newblock {Traveling waves for a thin film with gravity and insoluble
  surfactant}.
\newblock {\em SIAM J. Appl. Dyn. Syst.}, 14(4):1991--2012, 2015.

\bibitem{ELM11}
J.~Escher, Ph. Lauren\c{c}ot, and B.-V. Matioc.
\newblock Existence and stability of weak solutions for a degenerate parabolic
  system modelling two-phase flows in porous media.
\newblock {\em Ann. Inst. H. Poincar\'e Anal. Non Lin\'eaire}, 28(4):583--598,
  2011.

\bibitem{EL17}
J.~Escher and Ch. Lienstromberg.
\newblock {Travelling waves in dilatant non-Newtonian thin films}.
\newblock {\em J. Differential Equations}, 264(3):2113--2132, 2018.

\bibitem{escher2012modelling}
J.~Escher, A.-V. Matioc, and B.-V. Matioc.
\newblock Modelling and analysis of the {M}uskat problem for thin fluid layers.
\newblock {\em Journal of Mathematical Fluid Mechanics}, 14(2):267--277, 2012.

\bibitem{escher2013thin}
J.~Escher, A.-V. Matioc, and B.-V. Matioc.
\newblock Thin-film approximations of the two-phase {S}tokes problem.
\newblock {\em Nonlinear Analysis: Theory, Methods \& Applications}, 76:1--13,
  2013.

\bibitem{escher2013existence}
J.~Escher and B.-V. Matioc.
\newblock Existence and stability of solutions for a strongly coupled system
  modelling thin fluid films.
\newblock {\em NoDEA: Nonlinear Differential Equations and Applications}, pages
  1--17, 2013.

\bibitem{escher2014non}
J.~Escher and B.-V. Matioc.
\newblock Non-negative global weak solutions for a degenerated parabolic system
  approximating the two-phase {S}tokes problem.
\newblock {\em Journal of Differential Equations}, 256(8):2659--2676, 2014.

\bibitem{gan1}
F. Gancedo, E. Garcia-Juarez, N. Patel, and R. Strain.
\newblock {On the Muskat problem with viscosity jump: Global in time results}.
\newblock {\em Advances in Mathematics}, 345:552--597, 2019.

\bibitem{GW06}
H.~Garcke and S.~Wieland.
\newblock {Surfactant spreading on thin viscous films: nonnegative solutions of
  a coupled degenerate system}.
\newblock {\em SIAM J. Math. Anal.}, 37(6):2025--2048, 2006.

\bibitem{GG}
D.~P. Gaver and J.~B. Grotberg.
\newblock {The dynamics of a localized surfactant on a thin film}.
\newblock {\em J. Fluid Mech.}, 214:127--148, 1990.

\bibitem{granero2019global}
R. Granero-Belinch{\'o}n and M. Magliocca.
\newblock Global existence and decay to equilibrium for some crystal surface
  models.
\newblock {\em Discrete \& Continuous Dynamical Systems-A}, 39(4):2101--2131,
  2019.

\bibitem{granero2018asymptotic}
R. Granero-Belinch{\'o}n and S. Scrobogna.
\newblock On an asymptotic model for free boundary darcy flow in porous media.
\newblock {\em arXiv preprint arXiv:1810.11798}, 2018.

\bibitem{greenspan1978motion}
H. Greenspan.
\newblock On the motion of a small viscous droplet that wets a surface.
\newblock {\em Journal of Fluid Mechanics}, 84(1):125--143, 1978.

\bibitem{JG2}
J.~B. Grotberg and O.~E. Jensen.
\newblock {The spreading of heat or soluble surfactant along a thin liquid
  film}.
\newblock {\em Physics of Fluids A}, 5(58), 1993.

\bibitem{JG}
O.~E. Jensen and J.~B. Grotberg.
\newblock {Insoluble surfactant spreading on a thin viscous film: shock
  evolution and film rupture}.
\newblock {\em J. Fluid Mech.}, 240:259--288, 1992.

\bibitem{laurenccot2013gradient}
Ph. Lauren{\c{c}}ot and B.-V. Matioc.
\newblock A gradient flow approach to a thin film approximation of the {M}uskat
  problem.
\newblock {\em Calculus of Variations and Partial Differential Equations},
  47(1-2):319--341, 2013.

\bibitem{laurencot2014thin}
Ph. Lauren{\c{c}}ot and B.-V. Matioc.
\newblock A thin film approximation of the {M}uskat problem with gravity and
  capillary forces.
\newblock {\em Journal of the Mathematical Society of Japan}, 66(4):1043--1071,
  2014.

\bibitem{laurenccot2017finite}
Ph. Lauren{\c{c}}ot and B.-V. Matioc.
\newblock Finite speed of propagation and waiting time for a thin-film {M}uskat
  problem.
\newblock {\em Proceedings of the Royal Society of Edinburgh Section A:
  Mathematics}, 147(4):813--830, 2017.

\bibitem{laurencot2017self}
Ph. Lauren{\c{c}}ot and B.-V. Matioc.
\newblock Self-similarity in a thin film {M}uskat problem.
\newblock {\em SIAM Journal on Mathematical Analysis}, 49(4):2790--2842, 2017.

\bibitem{matioc2012non}
B.-V. Matioc.
\newblock Non-negative global weak solutions for a degenerate parabolic system
  modelling thin films driven by capillarity.
\newblock {\em Proceedings of the Royal Society of Edinburgh Section A:
  Mathematics}, 142(5):1071--1085, 2012.

\bibitem{pernas2019analysis}
T. Pernas-Casta{\~n}o and J. Vel{\'a}zquez.
\newblock Analysis of a thin film approximation for two-fluid Taylor-Couette
  flows.
\newblock {\em arXiv preprint arXiv:1905.13606}, 2019.

\bibitem{R1}
M.~Renardy.
\newblock {On an equation describing the spreading of surfactants on thin
  films}.
\newblock {\em Nonlinear Anal.}, 26(7):1207--1219, 1996.

\bibitem{sheludko1967thin}
A.~Sheludko.
\newblock Thin liquid films.
\newblock {\em Advances in Colloid and Interface Science}, 1(4):391--464, 1967.

\bibitem{simon1986compact}
J.~Simon.
\newblock {Compact sets in the space $L^p(0,T;B)$}.
\newblock {\em Ann. Mat. Pura Appl. (4)}, 146:65--96, 1987.

\bibitem{WCM}
M.~R.~E. Warner, R.~V. Craster, and O.~K. Matar.
\newblock {Fingering phenomena associated with insoluble surfactant spreading
  on thin liquid films}.
\newblock {\em J. Fluid Mech.}, 510:169--200, 2004.

\bibitem{HY}
S.~G. Yiantsios and B.~G. Higgins.
\newblock {A mechanism of Marangoni instability in evaporating thin liquid
  films due to soluble surfactant}.
\newblock {\em Physics of Fluids}, 22(022102), 2010.

\end{thebibliography}

\end{document}